\documentclass[12pt,reqno]{amsart}

\usepackage[dvips]{graphicx} 

\usepackage{%amssymb, latexsym, amsmath, amscd, array, 
hyperref
}

\usepackage{breakurl}   %this allows for url breaks in pdf

\theoremstyle{definition}

\numberwithin{equation}{section}

\newcommand\N {{\mathbb N}}

\title{Two-track depictions of Leibniz's fictions}

\author[M. Katz]{Mikhail G. Katz} \address{M.~Katz, Department of
  Mathematics, Bar Ilan University, Ramat Gan 5290002 Israel}
\email{katzmik@math.biu.ac.il}

\author[K. Kuhlemann]{Karl Kuhlemann}\address{K. Kuhlemann, Gottfried
  Wilhelm Leibniz University Hannover, D-30167 Hannover, Germany}
\email{kus.kuhlemann@t-online.de}

\author[D. Sherry]{David Sherry} \address{D. Sherry, Department of
  Philosophy, Northern Arizona University, Flagstaff, AZ 86011, US}
\email{David.Sherry@nau.edu}

\author[M. Ugaglia]{Monica Ugaglia} \address{M. Ugaglia, Il Gallo
  Silvestre, Localit\`a Collina 38, Montecassiano, Italy}
\email{monica.ugaglia@gmail.com}

\author[M. van Atten]{Mark van Atten} \address{M. van Atten, Paris
  Archives Husserl (CNRS / ENS), 45 rue d’Ulm, 75005 Paris, France}
\email{mvanatten@ens.fr}


\begin{document}

%\doublespacing
\thispagestyle{empty}

%\huge

\begin{abstract}
Leibniz described imaginary roots, negatives, and infinitesimals as
useful fictions.  But did he view such `impossible' numbers as
mathematical entities?  Alice and Bob take on the labyrinth of the
current Leibniz scholarship.
\end{abstract}


\maketitle

\tableofcontents

%\today

\section{What was Leibniz's take on `impossible' numbers?}

Leibniz described imaginary roots, negatives, and infinitesimals as
useful fictions.  But did he view such `impossible' numbers as
mathematical entities?  Did he envision a violation of the Archimedean
axiom?  And what were his `bounded infinities'?  Can a person of
infinite age have been born?  Did mathematical existence have
comparable meaning to Leibniz as to Hilbert?

Numbers that we take for granted today, such as negatives,
irrationals, imaginary roots, and infinitesimals, go beyond the
conceptual world of Greek mathematicians.%
\footnote{Some prominent 20th century scholars were still opposed to,
  or uncomfortable with, both irrationals and infinitesimals.  Thus,
  Errett Bishop opposed both the classical development of the real
  numbers and the use of infinitesimals in teaching calculus
  \cite{bish}.  For a discussion, see \cite{11a}, \cite{15a},
  \cite{18j}, \cite{17i}.}
In a sense such numbers are impossible, or fictional.

Many 17th century pioneers saw their task as either explaining, or
expanding upon, Greek mathematics.  A typical case is Fermat's
reconstruction of Apollonius' \emph{Plane Loci} (see \cite{Re94},
\cite{20e} for a discussion).  In particular, they had to justify the
status of certain `fictions' admitting no geometric representation.
Throughout his mathematical career, Leibniz argued for the virtues of
expanding the scope of `quantity' to include negatives, imaginaries,
and infinitesimals.  It is, however, not always easy to discern the
precise nature of Leibniz's attitude toward expanding the conceptual
resources of mathematics.

Leibniz saw Galileo's paradox of the infinite%
\footnote{Galileo observed that the natural numbers admit a one-to-one
  correspondence with the squares.}
as an indication that the concept of an \emph{infinite whole} is
contradictory because it contradicts the part-whole principle.%
\footnote{The principle asserts that a (proper) part is smaller than
  the whole.  The principle goes back to Euclid (\emph{Elements} I,
  common notion 8).}
Leibniz's view is in sharp contrast with the modern one, comfortable
with the concept of an infinite cardinality.  Leibniz often used
`infinite number' in the sense of what we would refer to today as
\emph{infinite cardinality} (contradictory in Leibniz's view),
indicating that he used the term `infinite number' in a generalized
sense.  When he spoke of the reciprocals of infinitesimals used in his
calculus, he tended to use either `infinite quantity' or
\emph{infinitum terminatum} (bounded infinity) rather than `infinite
number', though occasionally he used the latter term as well, as when
he defined an infinitesimal as an ``infinitely small fraction, or one
whose denominator is an infinite number.''%
\footnote{``fraction infiniment petite, ou dont le denominateur soit
  un nombre infini'' \cite[p.\;93]{Le02}.}

In his 1683 \emph{Elementa Nova matheseos universalis}, Leibniz
explained that some mathematical operations cannot be performed in
actuality, but nonetheless one can exhibit `a construction in our
characters' (\emph{in~nostris characteribus} \cite[p.\;520]{Le83}),
meaning that one can carry out a formal calculation, such as those
with imaginary roots.  Just as Leibniz is pushing the envelope by
extending the domain of quantities to include inassignable ones, he is
pushing the envelope by extending the meaning of `construction' to
include a mental operation using `our characters'.

Leibniz referred to infinitesimals as \emph{fictional entities}.  But
what is the precise meaning of that expression?  The crux of the
matter is whether Leibniz viewed infinitesimals as mathematical
entities.

Alice holds that the Leibnizian term `infinitesimal' does not refer to
a mathematical entity,%
\footnote{``[Numbers or ideal entities] are entities that are referred
  to.  Fictions, on the other hand, are not entities to which we
  refer.  \emph{They are not abstract entities}'' \cite[p.\,100]{Is90}
  (emphasis added).

``Leibniz conceived of infinitely small as compendia cogitandi for
  proofs and discovery and \emph{not as genuine mathematical
    entities}'' \cite[p.\;360]{Ra15} (emphasis added).  

``[R]eference to the infinite and infinitely small \emph{does not
    amount to the acceptance of genuine infinite entities}, but is a
  `way of speaking' referring ultimately to the only existing
  mathematical quantities, that is, finite quantities''
  \cite[p.\;441]{Ra20} (emphasis added).}
and sees Leibniz's expression \emph{fictional entities} as including
terms that only seem to refer to mathematical entities but in
actuality do not.  Alice reads the epithet \emph{fictional} as
undermining the noun \emph{entity}.  Furthermore, to Alice
`infinitesimals' do not refer to mathematical entities because such
would be inconsistent, i.e., contradictory, and more specifically
contrary to the part-whole axiom.

Bob holds that infinitesimals \emph{are} mathematical entities, and
interprets the expression \emph{fictional entities} as describing
entities of a special kind, namely \emph{fictional}.  Bob reads the
epithet \emph{fictional} not as undermining but as delimiting the
meaning of \emph{entity}, and views these mathematical entities as
consistent, as any mathematical entity would have to be; in
particular, they do not contradict the part-whole axiom.  Bob holds
that their fictionality references the fact that they are merely
accidentally impossible (in accordance with the Leibnizian philosophy
of knowledge) but nonetheless consistent, and therefore legitimate
mathematical entities in Leibniz's view.

\section
{Law and \emph{fictio juris}}

Did Leibniz view the term `infinitesimal' as tied up with
contradiction?  Alice cites as evidence the fact that Leibniz
sometimes used contradictory notions in jurisprudence.%
\footnote{``[E]ven though its concept [infinitesimal] may contain a
  contradiction, it can nevertheless be used to discover truths,
  provided a demonstration can (in principle) be given to show that
  its being used according to some definite rules will avoid
  contradiction.  This strategy of using `fictions' is not limited to
  mathematics and was very widespread in Law, the discipline which
  Leibniz first learned as a student, where it took the form of the
  `\emph{fictio juris}'{}'' (\cite[p.\;407]{Ra20}; emphasis in the
  original). }

Bob argues that jurisprudence fails to provide convincing evidence as
far as Leibniz's mathematical practice is concerned.  Bob holds that
non-contra\-diction was the very foundation of the mathematical method
for Leibniz%
\footnote{See \cite[\S\,3.4]{21a}.}
(see also Section~\ref{s7}), barring any inference from legal usage.

\section{Reference to violation of Euclid V.4}

Infinitesimals, as usually conceived, involve a violation of the
Archi\-me\-dean property.  One can therefore ask whether Leibniz ever
alluded to such a violation in writing.  In fact, Leibniz wrote in a
14/24 june 1695 letter to l'Hospital:
\begin{enumerate}\item[]
I use the term \emph{incomparable magnitudes} to refer to [magnitudes]
of which one multiplied by any finite number whatsoever, will be
unable to exceed the other, in the same way [adopted by] Euclid in the
fifth definition of the fifth book [of \emph{The
    Elements}].~\cite{Le95a}
\end{enumerate}
In modern editions of \emph{The Elements}, the notion of comparability
appears in Book~V, Definition~4.%
\footnote{``Magnitudes which when multiplied can exceed one another
  are said to have a ratio to one another" [Transl. Mueller].  A
  system of magnitudes satisfying Euclid\;V def.\;4 is said to be
  Archimedean, in connection with the so-called Archimedean axiom; in
  modern notation: ``for every $A$ and every $B$, a multiple $nA$ of
  $A$ exists such that $nA > B$".}
A similar discussion of incomparability in the context of Euclid's
definition appears in a 1695 publication of Leibniz's \cite{Le95b} in
response to Nieuwentijt's criticism.

Alice reads the Leibnizian reference to Euclid's Definition V.4, and
the violation thereof by infinitesimals when compared to ordinary
magnitudes, as merely a `nominal definition'.%
\footnote{``If one wants to infer existence, one cannot just rely on
  the nominal definition of `incomparables' (as not respecting the
  definition of Archimedean quantities)'' \cite[p.\;433]{Ra20}.}
Alice quotes Leibniz to the effect that nominal definitions could
harbor contradictions.  Alice holds that the true meaning of
infinitesimals resides in the Archimedean exhaustion-style unwrapping
of ostensibly infinitesimal arguments.%
\footnote{``The strict proof operating only with assignable quantities
  justifies proceeding by simply appealing to the fact that~$dv$ is
  incomparable with respect to~$v$: in keeping with the
  \emph{Archimedean axiom}, it can be made so small as to render any
  error in neglecting it smaller than any given''
  \cite[pp.\;567--568]{Ar13}.}

Bob argues that Archimedean paraphrases in exhaustion style constitute
an \emph{alternative} method rather than an unwrapping of the
infinitesimal method.%
\footnote{See \cite[\S\,1.3]{21a}.}
Bob notes that, while Leibniz warned that nominal definitions
\emph{may} harbor contradictions, there is no indication that they
\emph{must} do so; hence, regardless of whether one interprets the
violation of Euclid's Definition\;V.4 as a `nominal' move,
infinitesimals can still be consistent mathematical entities.

\section
{Fictions, useful fictions, and well-founded fictions}

Do fictions involve contradictions?  Some Leibnizian texts shed light
on the matter.  In 1674, Leibniz analyzed the area under the
hyperbola, and concluded that 
\begin{enumerate}\item[]
``the infinite is not a whole, but only a fiction, since otherwise the
part would be equal to the whole'' \cite{Le23}, A VII 3, 468; october
1674.
\end{enumerate}

Alice quotes this text as evidence that Leibniz uses the term
\emph{fiction} to refer to a contradictory infinite whole.%
\footnote{The Leibnizian passage is quoted as evidence in \cite{Ar13}
  and \cite{Ra20} as follows: 

``Even though this establishes the \emph{fictional nature of such
    infinite wholes}, however, this does not mean that one cannot
  calculate with them; only, the viability of the resulting
  calculation is contingent on the provision of a demonstration''
  (\cite[p.\;557]{Ar13}; emphasis added); 

``Here, the infinite area is that between the hyperbola and its
  asymptote (bounded on one side), and Leibniz argues that since
  taking it as a true whole leads to contradiction with the axiom that
  the whole is greater than its (proper) part, it should instead be
  regarded as a \emph{fiction}'' (\cite[p.\;405]{Ra20}; emphasis
  added).}

Bob points out that, while Leibniz does use the term \emph{fiction} in
this analysis of an infinite whole, he never refers to such
contradictory notions as either \emph{useful} or \emph{well-founded}
fictions; meanwhile, Leibniz does describe infinitesimals as both
useful fictions and well-founded fictions.%
\footnote{See \cite{13f}, \cite{14c}, \cite{18a}.}
Furthermore, Leibniz did not actually write that an infinite whole was
a fiction, contrary to Alice's inference.  Leibniz wrote that `the
infinite is not a whole, but only a fiction'.  That is not the same as
saying that an infinite whole is a fiction.  Therefore the inference
from the 1674 passage is inconclusive.

\section
{Infinite cardinalities and infinite quantities}
\label{s5}

Alice and Bob have argued about both the meaning of `infinite number'
in Leibniz, and his distinction between \emph{infinita terminata}
(bounded infinities) and \emph{infinita interminata} (unbounded
infinities).  One of the main sources for this Leibnizian distinction
is his \emph{De Quadratura Arithmetica}~\cite{Le04b}.

Leibniz's writings contain many speculations about the paradoxical
behavior of the \emph{infinita terminata}.  For example, Leibniz
mentioned the allegory of somebody of infinite age who nonetheless was
born; somebody who lives infinitely many years and yet dies
\cite[p.\;51]{Ar01}.  According to Leibniz, the kind of infinite
quantities one obtains by inverting infinitesimals is \emph{infinita
  terminata}, as in the example of an infinite-sided polygon.  Bob
argues that these ideas seem difficult today because of the prevalence
of a post-Weierstrassian mindset in traditional mathematical
training.%
\footnote{See Section~\ref{f11} for Skolem's formalisation of the idea
  of an infinite integer.}

Alice quotes passages where Leibniz argues that infinite wholes are
contradictory because contrary to the part-whole axiom.%
\footnote{``[Leibniz] argued in some critical comments on Galileo's
  \emph{Discorsi} in 1672 that the part-whole axiom must be upheld
  even in the infinite.  It follows that it is impossible to regard
  `all the numbers' and `all the square numbers' as true wholes, since
  then the latter would be a proper part of the former, and yet equal
  to it, yielding a contradiction'' \cite[pp.\;405--406]{Ra20}.}\,%
\footnote{Arguably, Leibniz in fact possessed the means to see that
  the part-whole axiom and the existence of infinite wholes are not
  incompatible \cite{Va11}.}
Alice holds that \emph{infinite number} necessarily means
\emph{infinite whole}, infinitesimals are their inverses, and
therefore all are contradictory.

Bob analyzes the Leibnizian distinction between bounded infinity and
unbounded infinity, and points out that the latter is akin to
cardinality.%
\footnote{See \cite[\S\,2.2]{21a}.}
The former are the inverses of infinitesimals, constitute a notion
distinct from cardinalities, and involve no contradiction.%
\footnote{``[U]nlike the infinite number or the number of all numbers,
  for Leibniz infinitary concepts do not imply any contradiction,
  although they may imply paradoxical consequences''
  \cite[Section~7]{Es21}.}
Bob holds that the expression \emph{infinite number} in Leibniz is
ambiguous and could refer either to cardinalities (contradicting the
part-whole axiom), or to (noncontradictory) \emph{infinita terminata}.
For a modern illustration of infinita terminata see Section~\ref{f11}.

\section
{Bounded infinities from Leibniz to Skolem}
\label{f11}

We provide a modern formalisation of Leibniz's \emph{infinita
  terminata} in terms of the extensions of~$\N$ developed by Skolem in
1933 \cite{Sk33}.  Such an extension, say~$M$, satisfies the axioms of
Peano Arithmetic (and in this sense is indistinguishable from~$\N$).
Yet~$M$ is a \emph{proper} extension, of which~$\N$ is an initial
segment.  Such models are sometimes referred to as nonstandard models
of arithmetic; see e.g., \cite{Ka91}.  Each element of the complement
\[
M\setminus\N
\]
is greater than each element of~$\N$ and in this sense can be said to
be infinite.

Notice that, depending on the background logical system, one can view
Skolem's extensions as either `potentially' or `actually' infinite (of
course in the former case neither~$\N$ nor~$M$ exists as a completed
whole).  The sense in which elements of~$M\setminus\N$ are `infinite'
is unrelated to the Aristotelian distinction.  An element of
$M\setminus\N$ provides a modern formalisation of the \emph{infinita
  terminata}.

\section
{Leibniz's rebuttal of Bernoulli's inference from series}
\label{s6}

In a 24 february/6 march 1699 letter to Bernoulli \cite{Le99}, Leibniz
noted that the infinity of terms in a geometric progression does not
prove the existence of infinitesimals:
\begin{enumerate}\item[]
You do not reply to the reason which I have proposed for the view
that, given infinitely many terms, it does not follow that there must
also be an infinitesimal term. This reason is that we can conceive an
infinite series consisting merely of finite terms or of terms ordered
in a decreasing geometric progression. I concede the infinite
plurality of terms, but this plurality itself does not constitute a
number or a single whole.  \cite[p.\;514]{Le89}
\end{enumerate}
Leibniz used the distinction between a plurality and an infinite whole
to refute Bernoulli's attempted inference from the existence of
infinite series to the existence of infinitesimals, and reiterated his
position against viewing an infinite plurality as a whole (see
Section~\ref{s5}).

Alice argues that Leibniz's exchange with Bernoulli about infinite
series shows that Leibniz viewed infinitesimals and infinite
quantities as contradictory.%
\footnote{``This remains Leibniz's position into his maturity and both
  arguments are to be found, for example, in the correspondence with
  Bernoulli in 1698 {\ldots} That is, he held that the part-whole
  axiom is constitutive of quantity, so that the concept of an
  infinite quantity, such as an infinite number or an infinite whole,
  involves a contradiction'' \cite[p.\;406]{Ra20}.}

Bob notes that Leibniz stresses the distinction between infinite
cardinality and infinite quantity (reciprocal of infinitesimals).  Bob
argues that the exchange with Bernoulli precisely refutes Alice's
attempt to blend infinite cardinality and infinite quantity so as to
deduce the inconsistency of infinitesimals.  Bob holds that Leibniz
didn't blend cardinality and quantity: only Alice did.  Leibniz, on
the contrary, emphasized the distinction in order to refute
Bernoulli's inference.  Bob holds that Leibniz's rebuttal of
Bernoulli's inference does a serviceable job of refuting Alice's
inference concerning a purported inconsistency of fictional entities,
as well.

\section{Mathematical possibility}
\label{s7}

Among Leibniz’s preparatory material for his \emph{Characteristica
  Universalis}, we find the following definiton of possible, dating
approximately from 1678:
\begin{enumerate}\item[]
A possible thing is that which does not imply a contradiction.%
\footnote{``Possibile est quod non implicat contradictionem'' A VI-2
  p. 495.  The definition is an additon made in 1678 to a text dating
  from 1671--1672 [ibid., p.\;487].}
\end{enumerate}
The same definition appears in many writings, as for instance the
24~february/6 march 1699 letter to Bernoulli analyzed in Section
\ref{s6}, where Leibniz wrote:
\begin{enumerate}\item[]
Possible things are those which do not imply a contradiction.%
\footnote{``Possibilia sunt quae non implicant contradictionem''
  \cite{Le99}.}
\end{enumerate}
If even in the broader framework of the \emph{Characteristica
  Universalis}, a thing is possible as soon as it causes no
contradiction, then certainly in the narrower mathematical context,
the absence of contradiction is sufficient to guarantee that the thing
is possible.  And in fact, that Leibniz meant the principle of
non-contradiction to apply to mathematics is evident from his second
letter to Clarke, from 1715:
\begin{enumerate}\item[]
The great foundation of mathematics is the \emph{principle of
  contradiction or identity}, that is, that a proposition cannot be
true and false at the same time, and that therefore A is A and cannot
be not A.  This single principle is sufficient to demonstrate every
part of arithmetic and geometry, that is, all mathematical principles.%
\footnote{``Le grand fondement des Mathematiques est le \emph{Principe
    de la Contradiction}, ou de l’\emph{Identit\'e}, c'est \`a dire,
  qu'une Enontiation ne sauroit etre vraye et fausse en m\^eme temps,
  et qu’ainsi $A$~est~$A$, et ne sauroit etre non\;$A$.  Et ce seul
  principe suffit pour demontrer toute l’Arithmetique et toute la
  Geometrie, c’est \`a dire tous les Principes Mathematiques''
  (Leibniz \cite{Le75}, 7:355–356).}
\cite[p.\;7]{Le00}
\end{enumerate}
In itself, the identification of possibility with the principle of
non-contradiction is not a novelty, as already in his \emph{Summa
  theologiae}, Thomas Aquinas clearly explained the major consequences
of this assumption:
\begin{enumerate}\item[]
But what implies contradiction is not submitted to divine omnipotence,
because it cannot bear the qualification of possible.%
\footnote{``Ea vero quae contradictionem implicant, sub divina
  omnipotentia non continentur, quia non possunt habere possibilium
  rationem'' (\emph{Summa theologiae}, I, q.\;25, a.\;3).  }
\end{enumerate}
But while in the Middle Ages possibility, and hence non-contradiction,
was deemed to be a necessary condition for the existence of an entity,
but not a sufficient one (not every possibility is actualized), Bob
argues that in Leibniz’s mathematics the condition is also sufficient:
mathematical existence is equivalent to mathematical possibility, and
the latter is wholly determined by a (global) principle of
non-contradiction.%
%
%\footnote{See also \cite[pp.\;47, 55]{Va15}.}
%
Of course, this is not the case in physics, so that Leibniz can
introduce the notion of accidental impossibilities, namely notions
that are possible – and hence they exist in mathematics – but not
necessarily instantiated in \emph{rerum natura}.  Accordingly, Leibniz
held a non-contradiction view of mathematical existence that can be
seen as an early antecedent of Hilbert's Formalism.%
\footnote{This observation was first made by Dietrich Mahnke, writing
  contemporaneously with the development of Hilbert's formalism.  See
  e.g., \cite[pp.\;284-287]{Ma27}.}
Bob argues that, to be usable in mathematics, a concept must first and
foremost be non-contradictory, and that Leibniz's letter undercuts
Alice’s claim that Leibniz viewed infinitesimals as contradictory.

\section{A-track and B-track}

Alice (A) and Bob (B) represent a pair of rival depictions in the
scholarly debate concerning the interpretation of Leibniz’s fictional
quantities such as infinitesimals and their reciprocals.  A similar
debate exists with regard to Cauchy's infinitesimal analysis
\cite{21f}.

On the A-track reading, these quantities, just like infinite wholes
violating the part-whole axiom, were contradictory concepts; the
expression \emph{fictional entities} describing them harbors a
contradiction.  Consequently, this reading denies that infinitesimals
were the very basis of the Calculus; formulations that use them were
merely figures of speech, abbreviating the Archimedean unwrappings
thereof.

On the B-track reading, what Leibniz viewed as contradictory were only
infinite wholes (involving a contradiction with the part-whole axiom),
but not infinite and infinitesimal quantities.  The latter were useful
and well-founded fictions involving a violation of the Archimedean
property.  Their legitimacy as mathematical entities stemmed from
their consistency, in an early anticipation of Hilbert's Formalism.

\end{document}